\documentclass[12pt]{iopart}
\usepackage{iopams}
\usepackage{color}
\date{Version Dowker-3a.tex - 8 November 2011}
\def\black{\color{black}}

\def\qedbox{\hbox{$\rlap{$\sqcap$}\sqcup$}}
\def\XIx\langle#1\rangle{h(#1)}

\newtheorem{theorem}{Theorem}

\newtheorem{lemma}{Lemma}

\newtheorem{remark}{Remark}

\def\Id{\hbox{Id}}\def\Tr{\hbox{\rm Tr}}

\makeatother

\begin{document}
\title[Gauss-Bonnet Theorem]{Heat trace asymptotics and the Gauss-Bonnet Theorem for general connections}
\author{C. G. Beneventano}
\address{Departamento de F\'{\i}sica, Universidad Nacional de La Plata\\
Instituto de F\'{\i}sica de La Plata, CONICET-Universidad Nacional de La Plata\\
C.C.67, 1900 La Plata, Argentina}
\ead{gabriela@fisica.unlp.edu.ar}
\author{P. Gilkey}
\address{Institute of Theoretical Science, University of Oregon, Eugene OR 97403 USA}
\ead{gilkey@uoregon.edu}
\author{K. Kirsten}
\begin{address}{Department of Mathematics, Baylor University \\ Waco, TX 76798, USA}\end{address}
\begin{ead}{Klaus\_Kirsten@baylor.edu}\end{ead}
\author{E. M. Santangelo}
\address{Departamento de F\'{\i}sica, Universidad Nacional de La Plata\\
Instituto de F\'{\i}sica de La Plata, CONICET-Universidad Nacional de La Plata\\
C.C.67, 1900 La Plata, Argentina}
\ead{mariel@fisica.unlp.edu.ar}
\begin{abstract} We examine the local super trace asymptotics for the de Rham complex defined by an arbitrary
super connection on the exterior algebra. We show, in contrast to the situation in which the connection in question
is the Levi-Civita connection, that these invariants are generically non-zero in positive degree and that the critical
term is not the Pfaffian.
\end{abstract}
\maketitle
\section{Introduction}
\subsection{The Chern-Gauss-Bonnet Theorem}
Throughout this paper, we will let $\mathcal{M}:=(M,g)$ be a compact smooth
Riemannian manifold without boundary of dimension
$m$. Let
$\chi(M)$ be the Euler characteristic of $M$. If $m$ is odd, then $\chi(M)$ is zero and if $m$ is even, then $\chi(M)$ is given by
integrating a suitable expression in the curvature tensor which can be described as follows. Let $\{e_i\}$ be a local orthonormal
frame for the tangent bundle and let $R_{ijkl}$ denote the components of the curvature tensor. We adopt the {\it Einstein}
convention and sum over repeated indices. With our sign convention, the components of the {\it Ricci tensor} are given by setting
$\rho_{ij}=R_{ikkj}$ and the {\it scalar curvature} is given by contracting the Ricci tensor and setting $\tau=\rho_{ii}$. If
$m=2\bar m$, then the {\it Pfaffian} (or {\it Euler form}) is given by setting:
\begin{equation}\label{eqn-1}
E_m:=\frac1{(-8\pi)^{\bar m}\bar m!}g(e^{i_1}\wedge...\wedge e^{i_m},e^{j_1}\wedge...\wedge
e^{j_m})R_{i_1i_2j_1j_2}...R_{i_{m-1}i_mj_{m-1}j_m}\,.\phantom{...}
\end{equation}
We have, for example,
$$
E_2=\frac\tau{4\pi}\quad\hbox{and}\quad E_4=\frac{\tau^2-4|\rho|^2+|R|^2}{32\pi^2}\,.
$$
One sets $E_m=0$ if $m$ is odd. Let $dx$ denote the Riemannian unit of volume on $M$:
$$dx=gdx^1...dx^m\quad\hbox{where}\quad g:=\sqrt{\det(g_{ij})}\,.$$ Chern \cite{Ch44} established the following
result which generalizes the classical 2-dimensional Gauss-Bonnet formula:
\begin{theorem}
If $(M,g)$ is a compact Riemannian manifold, then
$$\chi(M)=\int_ME_mdx\,.$$
\end{theorem}

\subsection{Heat Trace Asymptotics} Let $D$
be an operator of Laplace type on the space of smooth sections $C^\infty(V)$ to some vector bundle $V$ over $M$, i.e. locally
$D$ has the form:
\begin{equation}\label{eqn-2}
D=-\left\{g^{ij}\Id\partial_{x_i}\partial_{x_j}+A^i\partial_{x_i}+B\right\}\,.
\end{equation}
Here $A^i$ and $B$ are suitable locally defined endomorphisms of the vector bundle $V$. For example, the scalar
Laplacian $\Delta_{\mathcal{M}}^0$ is of this form since we may express:
$$\Delta_{\mathcal{M}}^0=-g^{-1}\partial_{x_i}gg^{ij}\partial_{x_j}\,.$$
Let
$e^{-tD}$ be the fundamental solution of the heat equation. If
$f$ is a smooth ``smearing" function, which is used to localize matters, then there is a complete asymptotic expansion as
$t\downarrow0$ of the form:
$$\Tr_{L^2}\left\{fe^{-tD}\right\}\sim\sum_{n=0}^\infty t^{(n-m)/2}\int_Ma_n(x,D)f(x)dx\,.$$
We set $f=1$ to define the global heat trace asymptotics
$$a_n(D):=\int_Ma_n(x,D)dx$$
and expand
\begin{equation}\label{eqn-3}
\Tr_{L^2}\left\{e^{-tD}\right\}\sim\sum_{n=0}^\infty t^{(n-m)/2}a_n(D)\,.
\end{equation}
The local invariants $a_n(x,D)$ vanish for $n$ odd and for $n$ even are given by a local formula which is homogeneous of order $n$
in the derivatives of the total symbol of $D$. The asymptotic series given in Equation~(\ref{eqn-3})  was first studied by S.
Minakshisundaram and {\AA}. Pleijel \cite{MP49} for the scalar Laplacian and the existence of the full asymptotic series in a very
general setting is due to Seeley
\cite{Se68b,Se69a}. This asymptotic formula can be generalized to the case of manifolds with boundary if suitable boundary
conditions are imposed; Kennedy, Critchley, and Dowker
\cite{KCD80} is a seminal work in this regard -- the invariants $a_n$ for $n$ odd appear for manifolds with boundary so this is a
convenient formalism.

\subsection{The Heat Equation and the Chern-Gauss-Bonnet Formula}
We restrict for the moment to the case in which $D=\Delta_{\mathcal{M}}^p$
is the Laplacian on the space $C^\infty(\Lambda^p(M))$ of smooth
$p$-forms. McKean and Singer \cite{MS67} observed that
\begin{equation}\label{eqn-4}
\sum_{p=0}^m(-1)^p\Tr_{L^2}(e^{-t\Delta_{\mathcal{M}}^p})=\chi(M)\,.
\end{equation}
They wrote: ``{\it Because of the complete cancellation of the time-dependent part of the
alternating sum Z, it is natural to hope that some fantastic cancellation will also take place {\bf in the small, i.e.} in the
alternating pole sum}". In other words, they conjectured that the following identity involving the super trace holds:
\begin{equation}\label{eqn-5}
a_n(x,d+\delta):=\sum_{p=0}^m(-1)^p a_n(x,\Delta_{\mathcal{M}}^p)=\left\{\begin{array}{lll}
E_m&\hbox{if}&n=m\\
0&\hbox{if}&n<m\end{array}\right\}\,.
\end{equation}
Here $d$ denotes exterior differentiation,
$\delta$ denotes the dual, interior multiplication, and
$$\Delta_{\cal M}^p:=d_{p-1}\delta_{p-1}+\delta_pd_p\,.$$ McKean and Singer
established Equation~(\ref{eqn-5}) in the special case that $m=2$; subsequently, Patodi \cite{P71} established this identity for
arbitrary
$m$ and thereby gave a heat equation proof of the Chern-Gauss-Bonnet theorem. By imposing absolute boundary conditions, manifolds
with boundary can be considered \cite{G75}.

\subsection{The de Rham Complex}\label{sect-1.4} Introduce a $\mathbb{Z}_2$ grading on the
exterior algebra by setting
$$\Lambda^e(M)=\oplus_p\Lambda^{2p}(M)\quad\hbox{and
}\quad\Lambda^o(M)=\oplus_p\Lambda^{2p+1}(M)\,.$$
We then have an elliptic complex of Dirac type
$$
(d+\delta)_{\mathcal{M}}^{e/o}:C^\infty(\Lambda^{e/o})\rightarrow C^\infty(\Lambda^{o/e})
$$
where the associated operators of Laplace type are given by:
\begin{eqnarray*}
&&\Delta_{\mathcal{M}}^{e}=(d+\delta)_{\mathcal{M}}^{o}(d+\delta)_{\mathcal{M}}^{e}=\oplus_p\Delta_{\mathcal{M}}^{2p},\\
&&\Delta_{\mathcal{M}}^{o}=(d+\delta)_{\mathcal{M}}^{e}(d+\delta)_{\mathcal{M}}^{o}=\oplus_p\Delta_{\mathcal{M}}^{2p+1}\,.
\end{eqnarray*}

We now review a bit of the classical tensor calculus. Let $\xi\in
T^*M$ be a cotangent vector. Let $\hbox{ext}(\xi)$ denote left exterior multiplication by $\xi$ and let
$\hbox{int}(\xi)$ be the dual, interior multiplication. We then have:
$$\hbox{ext}(\xi)\omega:=\xi\wedge\omega\quad\hbox{and}\quad
  g(\hbox{ext}(\xi)\omega,\tilde\omega)=g(\omega,\hbox{int}(\xi)\tilde\omega)\,.$$
If $\nabla^g$ denotes the Levi-Civita connection, then exterior differentiation $d$ and the adjoint interior
differentiation $\delta$ are given by the formulas:
$$
d\omega=\hbox{ext}(dx^i)\nabla_{\partial_{x_i}}^g\omega\quad\hbox{and}\quad
\delta\omega=-\hbox{int}(dx^i)\nabla_{\partial_{x_i}}^g\omega\,.
$$
We define $\gamma(\xi):=\hbox{ext}(\xi)-\hbox{int}(\xi)$ to give the exterior algebra the structure of a {\it
Clifford module}:
$$\gamma(\xi)\gamma(\eta)+\gamma(\eta)\gamma(\xi)=-2g(\xi,\eta)\Id\,.$$
Let $\nabla$ be a connection on the exterior algebra bundle $\Lambda(M)$. We suppose that $\nabla$ is a {\it super connection},
i.e. that
$\nabla=\nabla^e\oplus\nabla^o$ restricts to define connections on $\Lambda^e$ and $\Lambda^o$ separately, but aside from that we
impose no restictions on
$\nabla$. In physics, these connections play a role in the context of dark matter and supergravity; see, e.g., \cite{arno,duff,nieu}.
For such a connection, one defines:
\begin{eqnarray*}
&&(d+\delta)_\nabla^{e/o}:=\gamma(dx^i)\nabla_{\partial_{x_i}}^{e/o}:C^\infty(\Lambda^{e/o})\rightarrow
C^\infty(\Lambda^{o/e})\,,\\ &&\Delta_\nabla^{e/o}:=(d+\delta)_\nabla^{o/e}(d+\delta)_\nabla^{e/o}\,,\hbox{ and}\\
&&a_n(x,(d+\delta)_\nabla):=a_n(x,\Delta_\nabla^e)-a_n(x,\Delta_\nabla^o)\,.
\end{eqnarray*}
Equation~(\ref{eqn-5}) generalizes to this setting to become:
\begin{theorem}\label{thm-2}
Adopt the notation established above. Then
$$\Tr_{L^2}(e^{-t\Delta_\nabla^e})-\Tr_{L^2}(e^{-t\Delta_\nabla^o})=\chi(M)\,.$$
\end{theorem}
While this is perhaps not surprising, we shall give the proof in Section~\ref{sect-3} as the usual proofs in the literature (see,
for example, the discussion in \cite{G95}) assume the operators in question are self-adjoint and it is necessary to remove this
restriction.

If $D$ is an operator of Laplace type on a smooth vector bundle
$V$, then
$$a_0(x,D)=(4\pi)^{-m/2}\dim(V)\,.$$
Thus, in particular
$$a_0(x,(d+\delta)_\nabla)=(4\pi)^{-m/2}\left\{\dim(\Lambda^e)-\dim(\Lambda^o)\right\}=0\,.$$
The local ``fantastic cancellation" does not pertain in this setting for the higher order heat trace asymptotics. We shall
establish the following result in Section~\ref{sect-3}:
\begin{theorem}\label{thm-3} Let $\mathcal{M}=(\mathbb{T}^m,g_0)$ be the flat torus. Fix a point $x\in\mathbb{T}^m$. There exists a
super connection
$\nabla$ on $\Lambda^*(M)$ so that $(d+\delta)_\nabla$ is self-adjoint and so
that $a_{2n}(x,(d+\delta)_\nabla)\ne0$ for
$n\ge1$ for some point $x$ of $\mathbb{T}^m$.
\end{theorem}

\begin{remark}\rm We note that $E_m$ vanishes for the flat metric; thus, in particular, if $m$ is even, then
$a_m(x,(d+\delta)_\nabla)\ne E_m$.
\end{remark}

If $\nabla$ is a connection on $TM$, we can extend $\nabla$ to act on the cotangent bundle by requiring that the duality
equation is satisfied:
$$X\langle Y,\omega\rangle=\langle\nabla_XY,\omega\rangle+\langle X,\nabla_X\omega\rangle
\ \forall X,Y\in C^\infty(TM),\ \omega\in C^\infty(T^*M)
$$
where $\langle\cdot,\cdot\rangle$ denotes the natural pairing between $TM$ and $T^*M$. There is then a unique extension of $\nabla$
to a connection on the full exterior algebra which satisfies the Leibnitz rule:
$$
\nabla_X(\omega_1\wedge\omega_2)=(\nabla_X\omega_1)\wedge\omega_2+\omega_1\wedge(\nabla_X\omega_2)
\ \forall\ X\in C^\infty(TM),\  \omega_i\in C^\infty(\Lambda M)\,.
$$
Since such an induced connection restricts to a connection on each
$\Lambda^pM$ separately, necessarily any induced connection is a super connection. However, since an arbitrary super connection
need not obey the Leibnitz rule, not every super connection is induced from an underlying connection on the tangent bundle; in
particular, the connections used to establish Theorem~\ref{thm-3} will not satisfy the Leibnitz property.

In this setting, one has a weaker
version of the McKean-Singer local vanishing theorem:
\begin{theorem}\label{thm-4} Let $\nabla$ be a connection on $\Lambda M$ which is induced from an underlying connection on $TM$.
Then
$a_n(x,(d+\delta)_\nabla)=0$ for $3n<2m$.
\end{theorem}

\subsection{General elliptic complexes}
Instead of considering the {\it de Rham complex}
one can consider more generally a first order partial differential operator
\begin{equation}\label{eqn-6}
A:C^\infty(V_1)\rightarrow C^\infty(V_2)
\end{equation}
where $V_1$ and $V_2$ are smooth vector bundles over $M$.
We assume $V_1$ and $V_2$ are equipped with fiber metrics and form the formal adjoint
$A^*:C^\infty(V_2)\rightarrow C^\infty(V_1)$. We say that Equation~(\ref{eqn-6}) defines an {\it elliptic complex of Dirac type} if
the associated second order operators
$D_1:=A^*A$ and
$D_2:=AA^*$ are of Laplace type. We define $$a_n(x,A):=a_n(x,D_1)-a_n(x,D_2)\,.$$
Note that in the case of the de Rham complex, the associated Laplacians are given by $D_1=\Delta_{\mathcal{M}}^{e}$ and
$D_2=\Delta_{\mathcal{M}}^{o}$. Equation~(\ref{eqn-5}) generalizes to this setting to yield
\begin{equation}\label{eqn-7}
\int_Ma_n(x,A)dx=\left\{\begin{array}{lll}
\hbox{index}(A)&\hbox{if}&n=m\\
0&\hbox{if}&n\ne m\end{array}\right\}\,.
\end{equation}

Atiyah, Padodi, and Singer \cite{APS73} showed that if $A$ was the operator of the signature complex with coefficients in an
auxiliary bundle $W$, then $a_n(x,A)=0$ for $n<m$ while $a_m(x,A)$ gave the classical formula in terms of the Hirzebruch
$L$ polynomial and the Chern forms of the coefficient bundle $W$; a similar result pertained for the twisted spin complex. This
result then led to a heat equation proof of the Atiyah-Singer theorem in full generality; the Lefschetz formulas and other results
have also been established using these methods. We refer to
\cite{G95} for more details; the field is a vast one and it is not possible to give a full account in a short note such as this.

If $(M,g,J)$ is a K\"ahler manifold, the same approach led to a proof of the Riemann-Roch formula. However, if the metric in
question was not K\"ahler, the local invariants of the heat equation did not give rise to the classical formula; there were
divergence terms \cite{G73a,GNP97}. Rather than vanishing for $n<m$, in fact one could construct examples where
$a_n(x,\partial+\delta^{\prime\prime})$ was non-zero for certain values of $n<m$. But it remained an open problem to
construct index problems where the local index formula did not vanish for $a_2$; such examples are provided by Theorem~\ref{thm-3}.

The following result, which generalizes Equation~(\ref{eqn-7}) to the situation not involving adjoints, will be used in the
proof of Theorem~\ref{thm-2}:
\begin{theorem}\label{thm-5}
Let $A_\pm:C^\infty(V_\pm)\rightarrow C^\infty(V_\mp)$ be first order partial
differential operators over a compact Riemannian manifold without boundary such that $D_\pm:=A_\mp A_\pm$ are operators of
Laplace type on $C^\infty(V_\pm)$. Set $ a_n(x,A_+,A_-):=a_n(x,D_+)-a_n(x,D_-)$. Then
$$\int_Ma_n(x,A_+,A_-)dx=\left\{\begin{array}{lll}
\hbox{index}(A_+)&\hbox{if}&n=m\\
0&\hbox{if}&n\ne m\end{array}\right\}\,.
$$
\end{theorem}

\subsection{Outline of the paper} In Section~\ref{sect-2}, we shall summarize briefly the theory of compact operators on
Hilbert space that we shall need; we omit the proofs as they are standard and refer instead to \cite{C90} and \cite{Wi} for
further details. In Section~\ref{sect-3}, we will use these results to establish Theorem~\ref{thm-2} and Theorem~\ref{thm-5}. In
Section~\ref{sect-4}, we recall some formulas for the heat trace asymptotics which we use in Section~\ref{sect-5} to complete the
proof of Theorem~\ref{thm-3}. In Section~\ref{sect-6}, we prove Theorem~\ref{thm-4}.

\section{Spectral Theory of compact operators in Hilbert space}\label{sect-2}

\begin{lemma}\label{lem-1}
Let $\sigma(K)$ denote the spectrum of a compact operator $K$ on an infinite dimensional Hilbert space.
\begin{enumerate}
\item Let $0\ne\mu\in\sigma(K)$. Then:
\begin{enumerate}
\item $\mu$ is an eigenvalue of $K$.
\item There exists an integer $\nu_K(\mu)$ so that $\ker(K-\mu)^{\nu_K(\mu)}=\ker(K-\mu)^{\nu_K(\mu)+1}$.
\end{enumerate}
\item The eigenvalues can only accumulate at $0$.
\item $0\in\sigma(K)$ and $\sigma(K)$ is countable.
\end{enumerate}
\end{lemma}

Let $0\ne\mu\in\sigma(K)$ where $K$ is a compact operator on an infinite dimensional Hilbert space. By Lemma~\ref{lem-1},
the eigenvalues only accumulate at $0$. The {\it Riesz projection} is defined by setting
$$\pi_K(\mu):=\frac1{2\pi\sqrt{-1}}\int_\gamma(\Theta-K)^{-1}d\Theta ,$$
where $\gamma$ is any simple closed curve about $\mu$ in the complex plane which encloses no other eigenvalues of $K$.
Let $E_K(\mu):=\hbox{Range}\{\pi_K(\mu)\}$. We set $\pi_K(\mu)=0$ and $E_K(\mu)=\{0\}$ if
$\mu\notin\sigma(K)$.
\begin{lemma}\label{lem-2} Adopt the notation established above. Then
\begin{enumerate}
\item If $0\ne\mu\in\sigma(K)$, then:
\begin{enumerate}
\item $\pi_K(\mu)^2=\pi_K(\mu)$.
\item $E_K(\mu)$ is a finite dimensional space.
\item $E_K(\mu)=\ker(\mu-K)^{\nu_K(\mu)}$.
\end{enumerate}
\item If $0\ne\mu_i\in\sigma(K)$ with $\mu_1\ne\mu_2$, then $\pi_K(\mu_1)\pi_K(\mu_2)=0$.
\item Let $E_K(0):=\cap_{0\ne\mu\in\sigma(K)}\ker(\pi_K(\mu))$. Then
\begin{enumerate}
\item $KE_K(0)\subset E_K(0)$.
\item$\sigma(K|_{E_K(0)})=\{0\}$.
\end{enumerate}\end{enumerate}
\end{lemma}

\section{Heat trace formulas for the index}\label{sect-3}
\subsection{Spectral theory of operators of Laplace type}
 Let $D$ be an operator of Laplace type on the space of smooth sections to a vector bundle $V$ over a compact Riemannian manifold.
Standard elliptic theory, see for example the discussion in \cite{G95},  shows that the spectrum of $D$ (viewed as an unbounded
operator in Hilbert space) is contained in a cone about the positive real axis with arbitrarily small slope. Thus, in particular,
there exists
$\kappa>>0$ so that the operator
$(D+\kappa)$ is invertible; $K:=(D+\kappa)^{-1}$ is then a compact operator on $L^2(V)$ (viewed as a Banach space). We have
$$
   D\phi=\lambda\phi\Leftrightarrow
(D+\kappa)\phi=(\lambda+\kappa)\phi\Leftrightarrow(D+\kappa)^{-1}\phi=(\lambda+\kappa)^{-1}\phi\,.
$$
Thus the spectrum of $D$ is a countable set which only accumulates at $\infty$. Set  $\sigma(D)=\{\lambda_n\}$, set
$\nu_D(\lambda_n):=\nu_K((\lambda_n+\kappa)^{-1})$, and set
$$
E_D(\lambda_n):=\{\phi:(D-\lambda_n)^{\nu_D(\lambda_n)}\phi=0\}=E_K((\lambda_n+\kappa)^{-1})\,.
$$
Elliptic regularity then shows $E_D(\lambda_n)\subset C^\infty(V)$. Note that $0\in\sigma(K)$  plays no role and
the analysis of Section~\ref{sect-2} shows there is a direct sum decomposition (which is not an orthogonal direct sum
decomposition in general) of the form:
$$L^2(V)=\oplus_nE_D(\lambda_n)\,.$$
It is then immediate that:
$$
\Tr_{L^2}(e^{-tD})=\sum_n e^{-t\lambda_n}\dim(E_D(\lambda_n))\,.
$$

\subsection{Proof of Theorem~\ref{thm-5}} Let $A_\pm:C^\infty(V_\pm)\rightarrow C^\infty(V_\mp)$ be first order partial
differential operators over a compact Riemannian manifold without boundary such that $D_\pm:=A_\mp A_\pm$ are operators of
Laplace type on $C^\infty(V_\pm)$. Since
$ A_\pm D_\pm=D_\mp A_\pm$,
$$A_\pm:E_{D_\pm}(\lambda_n)\rightarrow E_{D_\mp}(\lambda_n)\,.$$
If $\lambda_n\ne0$, then $D_\pm=A_\mp A_\pm$ is invertible on $E_{D_\pm}(\lambda_n)$; consequently, $A_\pm$ is an isomorphism.
The decomposition
$$ L^2(V_\pm)=E_{D_\pm}(0)\oplus_{\lambda_n\ne0}E_{D_\pm}(\lambda_n)$$
shows that
$$\ker(A_\pm)=E_{D_\pm}(0)\quad\hbox{and}\quad\hbox{range}(A_\mp)=\oplus_{\lambda_n\ne0}E_{D_\pm}(\lambda_n)\,.$$
Consequently,
$$\hbox{index}(A_+)=\dim\{E_{D_+}(0)\}-\dim\{E_{D_-}(0)\}$$
and the usual cancellation argument shows that:
\medbreak\qquad
$\displaystyle\Tr_{L^2}(e^{-tD_+})-\Tr_{L^2}(e^{-tD_-})$
\medbreak\qquad\qquad
$=\displaystyle\sum_ne^{-t\lambda_n}
\left[\dim\left\{E_{D_+}(\lambda_n)\right\}-\dim\left\{E_{D_-}(\lambda_n)\right\}\right]$
\medbreak\qquad\qquad
$=\displaystyle\dim E_{D_+}(0)-\dim E_{D_-}(0)=\hbox{index}(A_+)$.\hfill\qedbox

\subsection{The proof of Theorem~\ref{thm-2}} The condition that $A_\mp A_\pm$ are of Laplace type is a condition only on the
leading order symbols of
$A_\pm$; it is unchanged by $0^{th}$ order perturbations of these operators. Furthermore, since $a_m(x,A_+,A_-)$ is given by a
local formula,
$\hbox{index}(A_+)$ varies continuously if we perturb the $0^{th}$ order symbol; as the index is integer valued it is therefore
unchanged by $0^{th}$ order perturbations. Consequently, if we are dealing with the generalized de Rham complex, the index must in
fact be $\chi(M)$ and Theorem~\ref{thm-2} follows.\hfill\qedbox

\section{Local formulas for the heat trace asymptotics}\label{sect-4}
Let $D$ be an operator of Laplace type on $C^\infty(V)$. We adopt the notation of Equation~(\ref{eqn-2}). There is a unique
connection
$\nabla$ on $V$ and a unique endomorphism $E$ of $V$ so that we may express $D$ using the Bochner formalism
$$Du=-\left\{g^{ij}u_{;ij}+Eu\right\} ,$$
where $u_{;ij}$ denotes second covariant differentiation of $u$ with respect to the connection $\nabla$.  Let
$\Gamma_{uv}{}^w$ be the Christoffel symbols of the Levi-Civita connection and let $\omega$ be the connection $1$-form of
$\nabla$. Using the notation established above, one then has \cite{G95} that:
\begin{equation}\label{eqn-8}
\begin{array}{l}
\omega_i=\textstyle\frac12g_{ij}\left\{ A^j+g^{lk}\Gamma_{lk}{}^j\Id\right\},\\
E=B-\displaystyle g^{ij}\left\{\partial_{x_i}\omega_{j}+\omega_{i}\omega_{j}
    -\omega_{k}\Gamma_{ij}{}^k\right\}\,.\vphantom{\vrule height 11pt}
\end{array}\end{equation}
Let $\Omega$ be the curvature of the connection $\nabla$. The following is well known -- see, for example, the discussion in
\cite{G04} and the references therein:
\begin{equation}\label{eqn-9}\begin{array}{l}
a_0(x,D)=(4\pi)^{-m/2}\Tr\{\Id\},\\
a_2(x,D)=(4\pi)^{-m/2}\textstyle\frac16\Tr\{6E+\tau\Id\},\vphantom{\vrule height 11pt}\\
a_4(x,D)=(4\pi)^{-m/2}\textstyle\frac1{360}\Tr\{
   60E_{;kk}+60\tau E+180E^{2}+12\tau_{;kk}\Id\vphantom{\vrule height 11pt}\\
\qquad+5\tau^{2}\Id-2|\rho|^{2}\Id+2|R|^{2}\Id+30\Omega_{ij}\Omega_{ij}\}.\vphantom{\vrule height 11pt}
\end{array}\end{equation}
Partial information is available for all the terms
\cite{BGO90}:
\begin{equation}\label{eqn-10}
a_{2n}(x,D)=\frac{(-1)^{n}n!}{(2n+1)!}\left\{-n
\Delta^{n-1}\tau\Tr(\hbox{Id})-(4n+2)\Delta^{n-1}\Tr(E)\right\}+...\,.
\end{equation}
We refer to
\cite{Gi75b,Ta73} for a discussion of $a_6$ and note that formulae for $a_8$ and $a_{10}$
are available in this setting \cite{AmBeOc89,Av90,Ven98}.

\section{Super trace asymptotics for the de Rham complex on a flat torus}\label{sect-5}
Let $\mathcal{M}=(\mathbb{T}^m,g)$ be the flat $m$-dimensional torus with the usual periodic parameters $(x^1,...,x^m)$ so that
$g(\partial_{x_i},\partial_{x_j})=\delta_{ij}$;
the analysis is simplified by taking a flat structure but a similar analysis pertains in the general setting. Let $V=\Lambda(M)$ be
the exterior algebra and let
$\Xi:=+\Id$ on $\Lambda^e$ and
$\Xi:=-\Id$ on
$\Lambda^o$ be the chirality operator which defines the super trace. Let $\gamma^i:=\hbox{ext}(dx^i)-\hbox{int}(dx^i)$. Let
$\theta$ be the connection
$1$-form of a connection
$\nabla$ on
$\Lambda(\mathbb{T}^m)$; $\nabla$ is a super connection if and only if $\theta\Xi=\Xi\theta$.  Set
$$A_\nabla:=\gamma^i(\partial_{x_i}+\theta_i)\,.$$

\begin{lemma}\label{lem-3}
If $\nabla$ is a super connection on $\Lambda(\mathbb{T}^m)$, then:
\begin{enumerate}
\item $A_\nabla:C^\infty(\Lambda^{e/o})\rightarrow C^\infty(\Lambda^{o/e})$ is an elliptic complex of Dirac type.
\item The operator $A$ is self-adjoint if and only if $\theta_i^*\gamma^i=-\gamma^i\theta_i$.
\item $\Tr\{E\Xi\}=-\Tr\{\gamma^i\gamma^j\theta_{j/i}\Xi\}$.
\end{enumerate}
\end{lemma}

\medbreak\noindent{\bf Proof.}
The first two assertions are immediate from the definitions that we have given. Since $g_{ij}=\delta_{ij}$, we
can raise and lower indices freely and use Equation~(\ref{eqn-8}) to compute:
\medbreak\quad
$D=\gamma^i(\partial_{x_i}+\theta_i)\gamma^j(\partial_{x_j}+\theta_j)$
\smallbreak\qquad
$=-\left\{\partial_{x_i}^2-(\gamma^i\theta_i\gamma^j+\gamma^j\gamma^i\theta_i)\partial_{x_j}
-(\gamma^i\theta_i\gamma^j\theta_j+\gamma^i\gamma^j\theta_{j/i})\right\}$,
\smallbreak\quad
$\omega_j=-\frac12\left(\gamma^i\theta_i\gamma^j+\gamma^j\gamma^i\theta_i\right)$,
\smallbreak\quad
$E=-(\gamma^i\theta_i\gamma^j\theta_j+\gamma^i\gamma^j\theta_{j/i})
+\frac12(\gamma^i\theta_{i/j}\gamma^j+\gamma^j\gamma^i\theta_{i/j})$
\smallbreak\qquad\quad
$-\frac14(\gamma^i\theta_i\gamma^j+\gamma^j\gamma^i\theta_i)(\gamma^k\theta_k\gamma^j+\gamma^j\gamma^k\theta_k)$.
\medbreak\noindent
We now take the super trace and show that most of the terms vanish. We use the identities $\Xi\theta_j=\theta_j\Xi$,
$\gamma^i\Xi=-\Xi\gamma^i$, and $\Tr(AB)=\Tr(BA)$ to compute:
\medbreak\quad
$-\Tr\{\gamma^i\theta_i\gamma^j\theta_j\Xi\}
  =-\Tr\{-\gamma^i\theta_i\Xi\gamma^j\theta_j\}
  =-\Tr\{-\gamma^j\theta_j\gamma^i\theta_i\Xi\}=0$,
\smallbreak\quad
$\frac12\Tr\{\gamma^i\theta_{i/j}\gamma^j\Xi+\gamma^j\gamma^i\theta_{i/j}\Xi\}
=\frac12\Tr\{\gamma^i\theta_{i/j}\gamma^j\Xi+\gamma^i\theta_{i/j}\Xi\gamma^j\}$
\smallbreak\qquad
$=\frac12\Tr\{\gamma^i\theta_{i/j}\gamma^j\Xi-\gamma^i\theta_{i/j}\gamma^j\Xi\}=0$.
\smallbreak\quad
$-\frac14\Tr\{\gamma^i\theta_i\gamma^j\gamma^k\theta_k\gamma^j\Xi
  +\gamma^i\theta_i\gamma^j\gamma^j\gamma^k\theta_k\Xi
  +\gamma^j\gamma^i\theta_i\gamma^k\theta_k\gamma^j\Xi
+\gamma^j\gamma^i\theta_i\gamma^j\gamma^k\theta_k\Xi\}$
\smallbreak\qquad
$=-\frac14\Tr\{\gamma^i\theta_i\gamma^j\gamma^k\theta_k\gamma^j\Xi
  +\gamma^i\theta_i\gamma^j\gamma^j\gamma^k\theta_k\Xi
  -\gamma^j\gamma^j\gamma^i\theta_i\gamma^k\theta_k\Xi
-\gamma^i\theta_i\gamma^j\gamma^k\theta_k\gamma^j\Xi\}$
\smallbreak\qquad
$=-\frac14\Tr\{-\gamma^i\theta_i\gamma^k\theta_k\Xi+\gamma^i\theta_i\gamma^k\theta_k\Xi\}=0$.\hfill\qedbox

We take $\theta_1:=f(x_2)\gamma^1\gamma^2\Xi$. Then $\theta_1\Xi=\Xi\theta_1$ so $\theta$ defines a super connection. Furthermore,
\begin{eqnarray*}
&&\phantom{......a}\theta_1^*=f(x_2)\Xi\gamma^2\gamma^1=-f(x_2)\gamma^1\gamma^2\Xi,\\
&&-\theta_1^*\gamma^1=f(x_2)\gamma^1\gamma^2\Xi\gamma^1=\gamma^1f(x_2)\gamma^1\gamma^2\Xi=\gamma^1\theta_1
\end{eqnarray*}
so $A_\nabla$ is self-adjoint. We have
$$\Tr\{\gamma^i\gamma^j\theta_{j/i}\Xi\}=\Tr\{\gamma^2\gamma^1\gamma^1\gamma^2\Xi\Xi
\partial_{x_2}f\}=\dim\{\Lambda(\mathbb{T}^m)\}\partial_{x_2}f\,.$$ We choose $f$ so all the derivatives of $f$ are non-zero
at some point
$x$ of $\mathbb{T}^m$. We then have $\{\Delta^{n-1}\Tr(E)\}(x)\ne0$ for $n\ge1$. Theorem~\ref{thm-3} now follows from
Equation~(\ref{eqn-9}) and from Equation~(\ref{eqn-10}).\hfill\qedbox

\section{The proof of Theorem~\ref{thm-4}}\label{sect-6}
Let $g$ be a Riemannian metric on $M$. Let $\nabla$ be a connection on $TM$. We expand
$$(\nabla_{\partial_{x_i}}-\nabla^g_{\partial_{x_i}})\partial_{x_j}=\theta_{ij}{}^k\partial_{x_k}\,.$$
We extend $\nabla$ to a connection on all of $\Lambda M$ using the Leibnitz rule as discussed previously.
Fix a point $x$ of
$M$. We normalize the choice of the coordinate system so that
$g(x)_{ij}=\delta_{ij}$ and so that the first derivatives of the metric $g_{ij/k}:=\partial_{x_k}g_{ij}$ vanish at $x$.
We introduce formal variables $g_{ij/\alpha}$ and $\theta_{ij}{}^k{}_{/\alpha}$ for the derivatives $\partial_x^\alpha$ of the
metric and of the normalized connection $1$-form. We let
$\mathcal{I}_m$ be the space of all polynomials which are
invariant, i.e. which are independent of the particular coordinate system chosen for evaluation.

There is a natural grading which is defined on $\mathcal{I}_m$. We set
$$
\hbox{deg}(g_{ij/\alpha})=|\alpha|\quad\hbox{and}\quad\hbox{deg}(\theta_{ij}{}^k{}_{/\alpha}):=|\alpha|\,.
$$
Let
$\mathcal{I}_{m,n}\subset\mathcal{I}_m$ consist of those invariants which are homogeneous of order $n$.  This grading can also be
expressed in a coordinate free fashion by noting that an invariant is homogeneous of order $n$ if and only if
$P(c^2g)=c^{-n}P(g)$.
In particular, there are no non-trivial invariants of odd order. Thus we may decompose
$$\mathcal{I}_m=\oplus_n\mathcal{I}_{m,n}\,.$$

The proof of Theorem~\ref{thm-4} will rest upon the following result:
\begin{lemma}
Adopt the notation established above:
\begin{enumerate}
\item $a_n(x,(d+\delta)_\nabla)$ defines an
element $\varepsilon_{m,n}\in\mathcal{I}_{m,n}$.
\item There is a natural restriction map $r:\mathcal{I}_{m,n}\rightarrow\mathcal{I}_{m-1,n}\rightarrow0$.
\item $r(\varepsilon_{m,n})=0$.
\item $r$ is injective from $\mathcal{I}_{m,n}$ to $\mathcal{I}_{m-1,n}$ if $3n<2m$.
\end{enumerate}
\end{lemma}

\medbreak\noindent{Proof.}
Assertion (i) follows from the Seeley calculus; $a_n$ is given by a local formula which is homogeneous of degree $n$ in the
derivatives of the total symbol of the operator. H. Weyl's theory of invariants of the orthogonal group \cite{W46} permits us to
express the spaces
$\mathcal{I}_{m,n}$ tensorially in terms of contractions of indices in polynomials of the covariant derivatives of the curvature tensor and of
$\theta$ where the sum ranges over indices $1\le i_\mu\le m$. If $P\in\mathcal{I}_{m,n}$, then $r(P)$ is defined by letting
the corresponding sum range over indices $1\le i_\mu\le m-1$; it now follows that $r$ is surjective.

It is worth illustrating this argument with an example. The scalar
curvature
$$\tau_m:=\sum_{i,j=1}^mR_{ijji}$$
defines an element of $\mathcal{I}_{m,2}$ for any $m$. The sequence $\{\tau_m\}$ satisfies $r(\tau_m)=\tau_{m-1}$ and for this
reason one usually does not subscript to explicitly demonstrate the dimension $m$.

There is an alternative description of the restriction map which is useful. Let $(N,g_N)$ be a Riemannian manifold of dimension
$m-1$ which is equipped with a connection $\nabla_N$ on the tangent bundle giving rise to a structure
$\theta_N=\{\theta_i\}_{i=1}^{m-1}$. Let $s$ be the usual periodic parameter on the circle. We let $M:=N\times S^1$,
$g_M=g_N+ds^2$, and $\theta_M=\{\theta_i\}_{i=1}^{m-1}$ define $\nabla_M$; the connection $\nabla_M$ on $T(N\times S^1)$ being
flat in the $S^1$ direction. Since the structures are flat in $S^1$, we have
$$r(P)(N,g_N,\theta_N)(x)=P(M,g_M,\theta_M)(x,s)\quad\hbox{for any}\quad s\in S^1\,.$$

Clifford multiplication $\gamma(ds):=\hbox{ext}(ds)-\hbox{int}(ds)$ defines an isomorphism from $\Lambda^e(M)$ to $\Lambda^o(M)$
which intertwines the associated Laplacians since the structures are induced from an underlying structure on $N$. Thus the super
trace of the heat asymptotics cancels and Assertion (iii) follows.

Let $P\in\mathcal{I}_{m,n}$. If $A$ is a monomial of order $n$, let $c(P,A)$ be
the coefficient of $A$ in $P$; we say $A$ {\it is a monomial of P} if $c(P,A)\ne0$. Express $A$ in the form:
$$A=g_{i_1j_1/\alpha_1}...g_{i_\ell
j_\ell/\alpha_\ell}\theta_{k_1l_1}{}^{n_1}{}_{/\beta_1}...\theta_{k_tl_t}{}^{n_t}{}_{/{\beta_t}}\,.$$
By hypothesis
\begin{equation}\label{eqn-11}
|\alpha_i|\ge2\quad\hbox{for}\quad 1\le i\le\ell
\end{equation}
since we have normalized the coordinate system so the first derivatives of the metric vanish. \black
 Furthermore, the order is given
by
\begin{equation}\label{eqn-12}
n=|\alpha_1|+...+|\alpha_\ell|+(|\beta_1|+1)+...+(|\beta_t|+1)\,.
\end{equation}
Let
$\hbox{deg}_a(A)$ be the number of times that an index $a$ appears in $A$;
\begin{eqnarray*}
\hbox{deg}_a(A)&=&\delta_{i_1,a}+\delta_{j_1,a}+\alpha_1(a)+...+\delta_{i_\ell,a}+\delta_{j_\ell,a}+\alpha_\ell(a)\\
&+&\delta_{k_1,a}+\delta_{l_1,a}+\delta_{n_1,a}+\beta_1(a)+...+\delta_{k_t,a}+\delta_{l_t,a}+\delta_{n_t,a}+\beta_t(a)\,.
\end{eqnarray*}\black
If $A$ is a monomial of $P$, we replace $x_a$ by $-x_a$ to see that
\begin{equation}\label{eqn-13}
\hbox{deg}_a(A)\equiv0\quad\hbox{mod}\quad2\,.
\end{equation}
If $r(P)=0$, then $\hbox{deg}_m(A)>0$ since the restriction map is defined by setting to zero those monomials which involve the
final index. Permuting the coordinate indices and applying Equation~(\ref{eqn-13}) then yields:
\begin{equation}\label{eqn-14}
\hbox{deg}_a(A)\ge2\quad\hbox{for}\quad 1\le a\le m\,.
\end{equation}
Suppose $r(P)=0$. Let $A$ be a monomial of $P$. We use Equation~(\ref{eqn-11}), Equation~(\ref{eqn-12}) and Equation~(\ref{eqn-14})
to estimate:
\begin{equation}\label{eqn-15}
\begin{array}{l}
2m\le\displaystyle\sum_{a=1}^m\hbox{deg}_a(A)=2\ell+3t+|\alpha_1|+...+|\alpha_\ell|+|\beta_1|+...+|\beta_t|\\
\phantom{aaaa}=2\ell+2t+n\le|\alpha_1|+...+|\alpha_\ell|+2t+n\le 2n+t\le 3n
\,.\vphantom{\vrule height 12pt}
\end{array}\end{equation}
This shows that $r$ is injective if $3n<2m$ which establishes Assertion (iv). \hfill\qedbox

\section*{Dedication} Stuart Dowker has made deep contributions to many areas of mathematics and physics and remains active
presently (see, for example, \cite{D10,D11}). The second and third authors have been honored to have collaborated with Stuart on
several papers (see, for example,
\cite{DGK99,DKG01}) and the remaining authors have used his results extensively. We
hope this paper serves as a fitting tribute to our colleague and friend and all the authors join in dedicating this paper to him.

\section*{Acknowledgements} Research of CGB and EMS was partially
supported by CONICET PIP1787, ANPCyT PICT909 and UNLP Proyecto 11/X492 (Argentina). Research of PG was partially supported by project MTM2009-07756
(Spain), by project 174012 (Serbia), and by DFG PI 158/4-6 (Germany). Research of KK was supported by the National Science
Foundation Grant PHY-0757791.

\end{document}